\documentstyle[12pt,epsf]{amsart}
\title[Explicit constructions of universal ${\Bbb R}$-trees]
{Explicit constructions of universal ${\Bbb R}$-trees and 
asymptotic geometry of hyperbolic spaces} 
\author{Anna Dyubina}
\address{School of Mathematical Sciences, Tel Aviv University, Ramat Aviv,
Israel}
\email{annadi@@math.tau.ac.il}
\author{Iosif Polterovich}
\address{Department of Mathematics, 
The Weizmann Institute of Science, Rehovot, Israel}
\email{iossif@@wisdom.weizmann.ac.il}
\def \phi{\varphi}
\def \epsilon{\varepsilon}
\numberwithin{equation}{subsection}
\theoremstyle{definition}
\newtheorem{definition}[equation]{Definition}           
\theoremstyle{plain}
\newtheorem{lemma}[equation]{Lemma} 
\newtheorem{theorem}[equation]{Theorem}

\begin{document}
\maketitle
\begin{abstract}
We present some explicit constructions of universal ${\Bbb R}$-trees
with applications to the asymptotic  geometry of hyperbolic spaces.
In particular, we show that any asymptotic cone  of a complete simply connected
manifold of negative 
curvature is a complete homogeneous ${\Bbb R}$-tree with the valency 
$2^{\aleph_0}$ at every point. It implies that all these asymptotic cones
are isometric depending neither on a manifold nor on an ultrafilter. 
It is also proved that the same tree can be isometrically embedded at 
infinity into such a manifold or into a non-abelian free group.
\end{abstract}
\section{Introduction and main results}
\subsection{Universal ${\Bbb R}$-trees as functional spaces}
The study of the ${\Bbb R}$-trees and their applications to different areas of
geometry and topology has been very intensive in the last years (see
a recent review [Bes] and references within). In particular, as it was
pointed out by Gromov, the ${\Bbb R}$-trees appear naturally in the asymptotic 
geometry of hyperbolic metric spaces ([Gr1],[Gr2]). The aim of this paper is
to present some explicit constructions related to this observation.

Recall that a metric space $T$ is {\it geodesic} (see [GH]) if
for every two points $t_1,t_2 \in T$ with the distance $a=d_T(t_1,t_2)$
between them, there exists an isometric inclusion $F:[0,a]\to T$,
such that $F(0)=t_1$, $F(a)=t_2$. The image $F([0,a])=[t_1,t_2]$ 
is called a {\it geodesic segment} joining the points $t_1$ and $t_2$.
 
A metric space $T$ is a {\it real tree} (or ${\Bbb R}$-{\it tree} for short) 
if any pair of its points can be connected with a unique geodesic segment,
and if the union of any two geodesic segments
$[x,y], [y,z]\subset T$ having the only
endpoint $y$ in common is the geodesic segment $[x,z]\in T$.

The {\it valency} of a point $t\in  T$ 
is the cardinal number of the set of connected
components of $T\setminus t$. An ${\Bbb R}$-tree is 
{\it $\mu$-universal} (see [MNO]) 
if any  ${\Bbb R}$-tree with the valency not greater than $\mu$ 
at every its point can be isometrically  embedded into it. 
We introduce an explicit construction of $\mu$-universal 
${\Bbb R}$-trees as spaces of functions.
 
For every cardinal $\mu \ge 2$ let $C_\mu$ be a set such that 
$\mbox{Card}(C_\mu)=\mu$ if $C_\mu$ is an infinite set and 
$\mbox{Card}(C_\mu)=\mu-1$ otherwise. Let us assume that each set
$C_\mu$ contains $0$ as one of it elements. 
\begin{definition}
\label{def1}
Let $A_\mu$ 
be a set of functions $f:(-\infty,\rho_f) \to C_\mu$ 
with the following properties:
1)~for every function $f$ there exists $\tau_f\le \rho_f$ such that 
$f(t)\equiv~0 \, \, \forall t<\tau_f$; 
2)~all functions $f$ are ``piecewise--constant'' 
from the right, i.e. for each $t\in (-\infty,\rho_f)$ there exists 
$\epsilon > 0$ such that 
$f|_{[t,t+\epsilon]}\equiv \mbox{const.}$
Denote the distance between two functions in $A_\mu$ by
\begin{equation}
\label{dA}
d_{A_\mu}(f_1,f_2)=(\rho_{f_1}-s)+(\rho_{f_2}-s),
\end{equation}  
where $s$ is the {\it separation moment} between the functions 
$f_1$,$f_2$ (see [PSh]):
$$
s=\sup \{ t|f_1(t')=f_2(t') \quad \forall t'<t \}. 
$$
\end{definition}
A relation similar to (\ref{dA}) in  [MNO] is called a 
"railroad-track" equation. 

For every $\mu$ the space $A_\mu$ is a 
complete real tree such that every its point has valency $\mu$ 
(see section 2.1). 
In particular, if $\mu=2$, i.e. if the set $C_\mu=\{0\}$,  the space $A_\mu$ 
is isometric to the real line. 
Functional spaces similar to the space $A_\mu$
were also considered in [PSh], [Ber], [Sh].
\begin{theorem} $\operatorname{(cf. [N], [MNO]).}$
\label{tree}
Consider any cardinal $\mu \ge 2$.

(i) (universality) The metric space $A_\mu$ is a $\mu$-universal
${\Bbb R}$-tree. 

(ii) (uniqueness) Any complete ${\Bbb R}$-tree with valency $\mu$ 
at every its point is isometric 

\hskip0.7truecm to  the space $A_\mu$.

(iii) (homogeneity) The metric space $A_\mu$ is homogeneous.
\end{theorem}
This theorem is proved in section 2.2 using simpler arguments 
than the analogous statements in ([N], [MNO]).

Universal
${\Bbb R}$-trees can be also constructed in a slightly different way. 
\begin{definition}
\label{def2}
Let $\mu$ be an infinite cardinal and let $\tilde A_\mu$ 
be a set of functions $f:[0,\rho_f) \to C_\mu$, $\rho_f\ge 0$ with the metric 
(\ref{dA}) which are piecewise-constant  from the right. In particular,
if $\rho_f=0$ we consider the function $f=f_{\emptyset}\in \tilde A_\mu$. 
For $\mu=2^{\aleph_0}$ we set 
$C_\mu=[0,1)$ and denote $\tilde A_\mu$ simply by $A$
(note that the definition of the space $A$ given in [DP] contains 
an inaccuracy).
\end{definition}
Such a representation of universal ${\Bbb R}$-trees
is more convenient for applications to asymptotic geometry of hyperbolic 
spaces.  Let us mention that for any infinite cardinal $\mu$ 
the spaces $A_\mu$ and $\tilde A_\mu$ are isometric due to part (ii)
of Theorem \ref{tree} being complete ${\Bbb R}$-trees with the valency
$\mu$ at every point. Note that this would fail for $\mu$ finite:
in this case the valency of the point $f_{\emptyset} \in \tilde A_\mu$ 
would be equal to  $\mu-1$.
\subsection{Isometric embedding at infinity} 
Let $X$ be a metric space with a distance function $d_X$. 
Intuitively, structure at infinity of the space $X$ is what is seen if one 
looks at the space $X$ from an infinitely far point (see [Gr2]).
One of the possible ways to treat this notion rigorously is as 
follows.
\begin{definition}
\label{embed}
A metric space $(T, d_T)$ admits an {\it isometric embedding at 
infinity} into the space $(X, d_X)$ if there exists
a sequence of positive scaling factors $\lambda_i \to \infty$ such that
for every point $t \in T$ one may find an infinite sequence 
$\{x_t^i\}$, $i=1,2,..$ of points in $X$
satisfying the relation
\begin{equation}
\lim_{i\to\infty}\frac{d_X(x_{t_1}^i, x_{t_2}^i)}{\lambda_i}=
d_T(t_1, t_2) 
\end{equation}  
for every $t_1,t_2 \in T$.
\end{definition}
In other words, every point of the space
$T$ corresponds to a sequence of points in $X$ tending 
to infinity, such that the ``normalized'' pairwise distances between the 
sequences in $X$ tend to the distances between the corresponding points 
in $T$.  Isometric embedding at infinity is a strengthening of the notion
of an {\it asymptotic subcone} of a metric space (see [Gr1], [GhH]).
\begin{theorem}
\label{main}
A $2^{\aleph_0}$-universal ${\Bbb R}$ tree can be isometrically embedded 
at infinity into:

(a) a complete simply connected manifold of negative curvature;

(b) a non-abelian free group.
\end{theorem} 
As usual, it is assumed in (a) that the curvature of a complete simply 
connected manifold is separated from zero. 
In the sequel we refer to such manifolds just as to 
{\it manifolds of negative curvature}. We prove part (a) in section 3.2.

Part (b) of this theorem is proved in section 3.1. We believe that its 
statement should remain true for any non-elementary hyperbolic group. 
Note that eventhough there exists a bilipschitz embedding of a non-abelian 
free group into a non-elementary hyperbolic group (see [Kap]),
this is not enough to generalize (b)
since isometric embeddings at infinity are not (at least, automatically)
preserved by bilipschitz mappings.

Let us mention that Theorem \ref{main} gives  a complete description of 
metric spaces admitting an isometric embedding at infinity into negatively 
curved manifolds or non-abelian free groups: any such (geodesic) space 
is a subtree of a $2^{\aleph_0}$-universal ${\Bbb R}$-tree 
(cf. [Gr1], [GhH]). 
\subsection{Asymptotic cones}
Let us recall the notion of an asymptotic cone of a metric space.
A {\it non-principal ultrafilter $\omega$ }
is a finitely additive measure on subsets of 
${\Bbb N}$ such that each subset has measure either $0$ or $1$ and all
finite subsets have measure $0$.
For any bounded function $h : {\Bbb N}\to {\Bbb R}$ its
limit $h(\omega)$ with respect to a non-principal ultrafilter  
$\omega$ is uniquely defined by the following condition: 
for every $\epsilon>0$ 
$$\omega(\{i\in I| \, |h(i)-h(\omega)\|<\epsilon\})=1.$$
\begin{definition} (See [Gr2], [DrW], [Dru]).
Fix a sequence of basepoints $O_i\in X$ and a sequence of scaling
factors $\lambda_i \to \infty$, $\lambda_i \in {\Bbb R}$.
Consider a set of sequences $g:{\Bbb N}\to X$ such that
$d_X(O_i, g(i))\le \mbox{const}_g \cdot \lambda_i$. 
To any pair of such sequences $g_1,g_2$ one may correspond a function
$$h_{g_1,g_2}(i)=\frac{d_X(g_1(i),g_2(i))}{\lambda_i}.$$
We say that the sequences $g_1$, $g_2$ are equivalent if the limit
$h_{g_1,g_2}(\omega)=0$. The set $T$ of all equivalence classes 
with the distance $d_T(g_1,g_2)=h_{g_1,g_2}(\omega)$ is
an {\it asymptotic cone} of $X$ with respect to the non-principal
ultrafilter $\omega$, sequence of basepoints $O_i$ and 
scaling factors $\lambda_i$. We denote it by 
$T=\mbox{Con}_{\omega}(X,O_i,\lambda_i)$.
\end{definition}
It was proved by M.~Gromov that any asymptotic cone of a hyperbolic
metric space is an ${\Bbb R}$-tree (see [Gr1],[Gr2])
\begin{theorem}
\label{con}
Any asymptotic cone of a manifold of negative curvature is isometric
to the $2^{\aleph_0}$-universal ${\Bbb R}$-tree.
\end{theorem}
We prove this theorem in section 2.3.
Theorem \ref{con} together with Theorem \ref{main} imply that any asymptotic 
cone of a manifold of negative curvature can be isometrically embedded at 
infinity into such manifold.

It follows from Theorem \ref{con} that 
asymptotic cones of manifolds of negative curvature depend neither 
on a manifold nor an ultrafilter 
or the sequence of scaling factors.
The statement of Theorem \ref{con} is folklore for
any non-elementary hyperbolic group (cf. [KapL], [Dru]). 
Let us note, however, that in general
asymptotic cones do depend on the choice of ultrafilters: in particular,
there exists a finitely generated group with 
non-homeomorphic asymptotic cones corresponding to different ultrafilters
(see [ThV]).

Main results of the present paper were announced in [DP].

\section{Universal ${\Bbb R}$-trees and asymptotic cones}
\subsection{Properties of the space $A_\mu$}
\begin{lemma}
\label{v}
The space $A_\mu$ is a real tree such that the valency of every its point 
is~$\mu$.
\end{lemma}
\noindent {\bf Proof.}
It is easy to see that the space $A_\mu$ with the metric (\ref{dA})
is a real tree.
Let us prove that the valency of each point is $\mu$.
Consider a function $f:(-\infty,\rho_f) \to C_\mu$  in $A_\mu$. For any 
$c \in C_\mu$ take a set
of functions $f_\delta:(-\infty,\rho_f+\delta) \to C_\mu$ such that
$f_\delta(t)= f(t)$ if $\quad -\infty < t < \rho_f$ and 
$f_\delta(t)\equiv c$ if $\quad \rho_f \le t < \rho_f+\delta.$
For different constants $c$ we get non--intersecting rays in $A_\mu$ starting 
from the point $f$ (note that if $\mu$ is finite the total number of rays
emanating from $f$ is exactly $\operatorname{Card}(C_\mu)+1=\mu$). 
Since $A_\mu$ is a real tree, these rays lie in different 
connected components of $A_\mu\setminus f$. 
Therefore, $A_\mu$ has valency $\mu$ at every its point. \qed
\begin{lemma}
\label{c}
The metric space $A_\mu$ is complete.
\end{lemma}
\noindent {\bf Proof.}
Consider a fundamental sequence $f_i:(-\infty,\rho_{f_i})\to C_\mu$ in $A_\mu$. 
Then $\rho_i=\rho_{f_i}$ is
also fundamental, since $d(f_i,f_j) \ge |\rho_i - \rho_j|$. Denote 
$\rho=\lim_{i\to \infty} \rho_i$.
Then 
$\lim_{i,j \to \infty}s_{ij}=\rho$, 
where $s_{ij}$ is a separation moment of functions $f_i,f_j$.
Indeed, $d_{A_\mu}(f_i, f_j)=(\rho_i - s_{ij})+(\rho_j - s_{ij})$ and hence
$$
s_{ij}=\frac{\rho_i + \rho_j - d_{A_\mu}(f_i,f_j)}{2} \to \rho.
$$
Therefore for each $\rho' < \rho$ the functions $f_i$ coincide for sufficiently
large $i$ on 
the interval $(-\infty,\rho')$. 
Define $f$ by $f(x)=\lim_{i\to \infty} f_i(x)$ for any
$x \in (-\infty,\rho)$. Note that $f \in A_\mu$, since for any 
$x \in (-\infty,\rho)$
there exist $I\in {\Bbb N}$ and $\epsilon >0$, such that 
$f|_{(-\infty,x+\epsilon]}\equiv f_I|_{(-\infty,x+\epsilon]}$ and hence $f$ is 
``piecewise--constant'' from the right.
\qed
\subsection{Proof of Theorem \ref{tree}}
(i) Let $\Psi$ be a tree of cardinality $card(\Psi)\le \mu$. 
We say that a subtree $T\subset \Psi$ is {\it good} if 
it satisfies the following two properties: a) $T$ is nonempty and connected;
b)$\forall t \in T$ either $T$ has a nonempty intersection with all connected 
components of the  set $\Psi\setminus t$ (in this case we say that $t$ is 
{\it full-valented}), or $T$ has a non--empty intersection with at most $2$
connected components of $\Psi\setminus t$ (then we call $t$ 
{\it few--valented}).
Note that union of the nested good subtrees is also a  good subtree.
Consider pairs $P=(T,F)$ such that $T$ is a good subtree and 
$F:T \to A_\mu$ is an isometric inclusion. Let us introduce a partial order on
such pairs: $(T_1,F_1) \le  (T_2,F_2)$ if $T_1 \subset T_2$ and 
$F_2|_{T_1}=F_1$. The set of such pairs is nonempty since it contains all
inclusions of a one-point set. 
Any linearly ordered set of pairs $P_\alpha=(T_\alpha,F_\alpha)$ has
a maximal element $(T',F')$,  where $T'=\cup T_\alpha$ and 
$F'|_{T_{\alpha}}=F_{\alpha}$. Then by Zorn's lemma there exists a maximal
element in the whole set of pairs: $P_{max}=(T_{max}, F_{max})$.
Let us show that $T_{max}=\Psi$.
If not, there exists $\gamma \in \Psi\setminus T_{max}$.
Take an arbitrary $t \in T_{max}$ and consider the geodesic segment
$[t,\gamma]$. The intersection $[t,\gamma]\cap T_{max}$ is either  a
closed segment $[t,t']$ or a semi-interval $[t,t')$.
Consider each of these cases separately.

{\it Case 1.} $T_{max}\cap [t,\gamma]=[t,t').$  Take $T_{max}'=T_{max}\cup t'$.
It is a good subtree since $T_{max}'$ intersects just with a single
connected component of $\Psi\setminus t'$ (since $T_{max}$ is connected),
and for the points from $T_{max}$ the condition b) remains true.
Due to completeness of $A_\mu$, the isometric inclusion of $T_{max}$ can be
continued to the isometric inclusion of $T_{max}'$. This contradicts with
the maximality of $P_{max}$. 

{\it Case 2.} $T_{max} \cap [t,\gamma]=[t,t']$. Then $T_{max}$ does not 
intersect
with the connected component of $\Psi\setminus t'$ containing $\gamma$,
therefore $t'$ is few--valented.
In each component of $\Psi\setminus t'$ non-intersecting with
$T_{max}$ we fix a point $t_\nu$ and consider geodesic segments $[t', t_\nu]$.
Take $T_{max}'=\cup_{\nu} [t',t_\nu]\cup T_{max}$.
This is a good subtree. Indeed, $t'$ is full-valented. 
For the points $T_{max}\setminus t'$ the condition b) remains true,
and the points of $T_{max}'\setminus T_{max}=\cup_\nu (t',t_\nu]$
are few--valented. 

  We build the continuation of $F$ to $T'_{max}$. Consider $F(t') \in A_\mu$.
The point $F(t')$ is an endpoint for $\mu$ rays lying in different
connected components of $A_\mu\setminus F(t')$.
Points of $F(T_{max})$ may lie in at most two of them.
To every segment $[t',t_\nu]$ we isometrically correspond a segment 
$[F(t'), F(t_\nu)]$ on one of the ``free'' rays (i.e. having
no points of $T_{max}$) of $A_\mu\setminus F(t')$. This is 
possible since cardinality of the set of the ``free'' rays is always
greater or equal than of the set of points $t_\nu$.
Once again, we get a contradiction with the maximality of $P_{max}$.

Therefore, $T_{max}=\Psi$ and $F_{max}$ is its isometric emebedding into
the space $A_\mu$. This completes the proof of (i).

\noindent (ii) Let $\Psi$ be a complete real tree such that the valency
of every its point has cardinality $\mu$.
We call the pair $(T,F)$ {\it good}, if $T$ is a good subtree and
$F$ preserves full-valentness. The same as in (i), the set of good pairs
is non-empty and one may introduce partial order on it. It may be shown
analogously that every linear ordered set of good pairs has a maximal element
and hence by Zorn's lemma there exists a maximal pair $(T_{max}, F_{max})$,
where $T_{max}=\Psi$. Let us prove that $F_{max}$ is a surjection.
If not, there exists $a \in A_\mu$ not lying in $F_{max}(T_{max})$.
Take an arbitrary point $p \in F_{max}(T_{max})$ and consider $[a,p] \cap
F_{max}(T_{max})$. Since $F_{max}(T_{max})$ is complete 
(being an isometric image of
a complete metric space $T_{max}=\Psi$), it is a closed segment
$[a,p']$. But in this case $F_{max}(T_{max})$ intersects not with all
connected components of 
$A_\mu\setminus p'$ while $F^{-1}_{max}(p') \in \Psi$ is clearly
full--valented --- we got a contradiction with the fact that $(T_{max},
F_{max})$ is a good pair.
Therefore, $F_{max}$ is an isometry between $A_\mu$ and~$\Psi$.

\medskip

\noindent (iii)
Let us note that in the proof
of (ii) we could consider only such good pairs $(T,F)$ 
that satisfy $(T_0, F_0)\le (T,F)$, 
where $T_0={\gamma_0}$, $F_0(\gamma_0)=a_0$
for some arbitrary fixed $\gamma_0 \in \Psi$ and $a_0 \in A_\mu$. 
Then finally we would obtain an
isometry between $\Psi$ and $A_\mu$ which moves $\gamma_0$ into $a_0$.
Setting $\Psi=A_\mu$ implies that $A_\mu$ is homogeneous. \qed

\subsection {Asymptotic cones of manifolds of negative curvature}
In this section we prove Theorem \ref{con} stating that any
asymptotic cone of a manifold of negative curvature is isometric
to the space $A$ (see Definition \ref{def2}).

\noindent {\bf Proof of Theorem \ref{con}.}
Let $O_i\in X$ be some sequence of basepoints and $\lambda_i\to \infty$ be
some sequence of scaling factors.
An asymptotic cone $\operatorname{Con}_\omega (X,O_i,\lambda_i)$
of a manifold of negative curvature $X$ 
is a complete geodesic metric  space.
Let us show that the valency at every its point is 
$2^{\aleph_0}$.

Consider a point $\xi \in \operatorname{Con}_\omega (X,O_i,\lambda_i)$ 
and let $\{x_i\}$ be 
one of the sequences in $X$ corresponding to this point.
For every $\alpha \in [0,2\pi)$ and for every $\rho>0$ consider a sequence
$\{x_i^{\alpha,\rho}\}$ constructed according to the following rule.
Consider the geodesic segment $[O_i,x_i]$.
The segment $[x_i, x_i^{\alpha,\rho}]$ is of length $\lambda_i\rho$ and the 
angle between it and $[O_i,x_i]$ is equal to $\alpha$.
It clearly follows from the triangle inequality that the distances
between $O_i$ and $x_i^{\alpha,\rho}$ grow not faster than linearly.
Let us show that the distance between two points of the asymptotic 
cone $\xi^{\alpha_j,\rho_j}$ 
corresponding to the sequences
$\{x_i^{\alpha_j,\rho_j}\}$, $j=1,2$, is equal to $\rho_1+\rho_2$
if $\alpha_1 \ne \alpha_2$ and $|\rho_1-\rho_2|$ if $\alpha_1=\alpha_2$
The latter is due to the fact that
$\alpha_1=\alpha_2$ implies 
$d_X(x_i^{\alpha_1,\rho_1},x_i^{\alpha_2,\rho_2})=\lambda_i|\rho_1-\rho_2|$.
If $\alpha_1 \ne \alpha_2$ we need the following lemmas.
\begin{lemma}
\label{26}
Let $ABC$ and $A'B'C'$ be geodesic triangles in $\Bbb H$ and in $X$ 
respectively, such that $|AB|=|A'B'|=u$, $|AC|=|A'C'|=v$, 
$|BC|=w$, $|B'C'|=w'$, $\angle A = \angle A'=\phi$, 
where $|\cdot|$ denotes the
length of a segment in the corresponding metric. 
Then  $w <w'$.
\end{lemma}
\noindent {\bf Proof.}
Consider a comparison triangle $PQR \in {\Bbb H}$ for $ABC$:
$|PQ|=u$, $|PR|=v$, $|QR|=w$ and let $\angle P=\psi$.
By Alexandrov comparison theorem (see  [GhH) $\psi>\phi$.
On the other hand by hyperbolic cosine rule (see [Bea])
$$\cosh w =\cosh u\cosh v-\cos \phi \sinh u \sinh v <$$
$$<\cosh u\cosh v-\cos \psi \sinh u \sinh v = \cosh w',$$
and hence $w<w'$. \qed
\begin{lemma}
\label{nep}
Consider a point $O\in X$ and two geodesic rays $l_1$ and $l_2$
starting from $O$ with a non-zero angle between them. Let $A_i\in l_1$
and $B_i\in l_2$ be two sequences of points such that 
$d_X(O,A_i)=\lambda_i\rho_1$ and $d_X(O,B_i)=\lambda_i\rho_2$,
where $\lambda_i \to \infty$.
Then 
\begin{equation}
\label{pum}
\lim_{i\to \infty}\frac{d_X(A_i,B_i)}{\lambda_i}=\rho_1+\rho_2.
\end{equation}  
\end{lemma}
\noindent {\bf Proof.}
Let  $O'A_i'B_i'$ be the triangle in ${\Bbb H}$ corresponding to
$OA_iB_i$ as in lemma \ref{26}.
Then due to triangle inequality and lemma \ref{26} we have:
$$
\rho_1+\rho_2 \ge \frac{d_X(A_i,B_i)}{\lambda_i} 
\ge \frac{d_{\Bbb H}(A_i',B_i')}{\lambda_i} \to
\rho_1+\rho_2,
$$
where the limit at the right follows from a simple asymptotic analysis of 
distance formulas on the hyperbolic plane.
This proves the relation (\ref{pum}). 
\qed

\smallskip

Lemma \ref{nep} implies that
\begin{equation}
\frac{d_X(x_i^{\alpha_1,\rho_1},x_i^{\alpha_2,\rho_2})}{\lambda_i} \to 
\rho_1+\rho_2
\end{equation}
Therefore $\{x_i^{\alpha_1,\rho_1}\}=\{x_i^{\alpha_2,\rho_2}\}$ iff
$\alpha_1=\alpha_2$ and $\rho_1=\rho_2$. For every $\alpha \in [0,2\pi)$
we get an infinite ray starting from the point $\xi$ and for different 
$\alpha$ these rays have the only point $\xi$ in common and hence lie
in different connected components of 
$\operatorname{Con}_{\omega}(X,O_i,\lambda_i)\setminus \xi$, i.e.
the valency of the point $\xi$ is $2^{\aleph_0}$.
Since $\xi$ is an arbitrary point in the asymptotic cone
$\operatorname{Con}_{\omega}(X,O_i,\lambda_i)$ 
we have proved that it is isometric to the metric space~$A$.~\qed

\section{Proof of Theorem \ref{main}}
\subsection{Part(b)}
We start with the easier part of Theorem \ref{main}, namely, we prove that 
the space $A$ can be isometrically embedded at infinity into a non-abelian
free group.

Let us fix  some sequence of  scaling factors $\lambda_i\to \infty$ 
(for simplicity we set $\lambda_i=i$).

Consider a Cayley graph corresponding to a subgroup generated
by two arbitrary generators $\gamma_1,\gamma_2$ of the group $\Gamma$,
and let $G$ be its part corresponding to non-negative powers of 
$\gamma_1,\gamma_2$.
Consider the standard coding of vertices of $G$ with binary sequences
of  $\gamma_1$-s and $\gamma_2$-s
taking the unity as the initial vertex. 
Setting $\gamma_1=0$ and $\gamma_2=1$ we correspond
to every vertex of $G$ a binary number in a unique way.

For any point $t\in [0,\rho_f)$  let $[t_l, t_r)\ni t$ be the maximal 
interval where the function $f\in A$ is constant. For every $\epsilon > 0$ 
consider the following ``cut-off''
function: $f^\epsilon(t)=f(t)$ if $t_r-t_l\ge \epsilon$ and 
$f^\epsilon(t)=0$ otherwise.

Fix some bijection $F$ from $[0,1)$ to the set of all binary sequences. 
Our aim is to build a sequence of vertices $\{ x_i \} \in G$ 
(which are actually elements of the group $\Gamma$) 
corresponding to an arbitrary function $f\in A$.
Consider a function 
$f^{1/\sqrt{i}}$, 
and let  
$l^i_1,...,l^i_{k_i}$ be the lengths of the segments on which it is constant,
and $0=a^i_1,...,a^i_{k_i}$ be the values of $f^{1/\sqrt{i}}$ on them.
Consider the binary sequence $F(a_1^i)$ and take its first $[l_1\cdot i]$ 
elements, where $[\cdot]$
denotes the integer part.
We get a binary number which corresponds to a unique point $x_1^i \in G$.
From the point $x_1^i$ we go along the sequence $F(a_2^i)$ for the time
$[l_2\cdot i]$ and so on until $x^i_{k_i}=x_i$.
Let us check that 
$$\frac{d_{G}(x_i,e)}{i} \to \rho_f.$$
 Note that $k_i \le 2\sqrt{i}\rho_f+1$ since the number of intervals
on which the function $f^{1/\sqrt{i}}$ does not vanish is not greater than
$\sqrt{i}\cdot \rho_f$.
Therefore 
$$d_{G}(x_i, e)=\sum_{r=1}^{k_i} [l^i_r\cdot i]\ge 
\sum_{r=1}^{k_i}(l^i_r \cdot i -1) = \rho_f\cdot i - k_i.$$
Thus, 
$$\rho_f =\sum_{r=1}^{k_i} l^i_r\ge 
\frac{d_{G}(x_i,e)}{i}\ge \rho_f-\frac{k_i}{i},$$
and since $k_i/i\to 0$ we get $d_{G}(x_i,e)/i \to \rho_f$.

Consider now two arbitrary functions $f,g \in A$. Let $s$ be their
separation moment, and let $s_i$ be moments of separation of the
segments $[e,x^f_i]$ and $[e,x^g_i]$.
Note that for sufficiently small $\epsilon$ the moment of separation
between $f^\epsilon$ and $g^\epsilon$ coincides with the separation moment 
of $f$ and $g$, and $f^\epsilon|_{[0,s]}=(f|_{[0,s]})^\epsilon.$
Applying previous arguments to the function $f|_{[0,s]}$ we get that
$d_{G}(s_i,e)/i \to s$, since for the sufficiently large $i$ we have
$s_i=x_i^{f|_{[0,s]}}.$ 
Hence $$\frac{d_{G}(x_i^f,x_i^g)}{i}=
\frac{d_{G}(x_i^f,e)}{i}+
\frac{d_{G}(x_i^g,e)}{i}-
\frac{d_{G}(s_i,e)}{i} \to \rho_f+\rho_g-2s=d_A(f,g).$$
This completes the proof of part (b) of the theorem.\qed
\subsection{Part(a)}
In this section we prove that the metric space $A$ can be isometrically 
embedded at infinity into any manifold of negative curvature.
First  let us mention that in case of the hyperbolic plane ${\Bbb H}$
(or, more general, in case of a complete simply connected manifold 
of {\it constant} negative curvature) 
one may isometrically embed the space $A$ into it in a very explicit and 
simple way. Consider the Poincare disc model of  ${\Bbb H}$ and let us
introduce the hyperbolic polar coordinates $(r, \phi)$ centered at
some $O\in {\Bbb H}$.
To every function $f(t)\ne f_\emptyset \in A$, $f:[0,\rho) \to [0,1)$  we
correspond a sequence $x_n=(\rho_i, \phi_i)\in {\Bbb H}$, where 
$\rho_i=\rho i$ and 
$$\phi_i=2\pi f(0)+\int_0^{\rho}e^{-ti}f(t)dt.$$ 
To the ``empty function'' $f_\emptyset$ we correspond a constant sequence
$x_i=O$. We also set $\lambda_i=i$.
Using some asymptotic analysis of distance formulas on the hyperbolic plane
one can prove that this is indeed an isometric embedding at infinity.

\smallskip

Let us prove the theorem in the general case. 
We divide the proof into several steps.
\subsection*{Step 1.}
Let us fix a sequence of  
scaling factors (as before, let $\lambda_i=i$). 
Let $O \in X$ be a basepoint and $\l_0 \in X$ some fixed geodesic ray
emanating from $O$. For every $0 \le \phi < 2\pi$ let us fix a 
ray $l_\phi \subset X$ starting from the point $O$ such that the angle
between $\angle(l_0,l_\phi)=\phi$. 
As in the case of the  hyperbolic plane, let us correspond to 
$f_\emptyset$ the constant sequence $x_i=O$.

We prove that for any $f:[0,\rho)\to [0,1)$ and for any $i>0$ 
there exists a unique naturally--parametrized path 
$\gamma_i: [0,\rho\cdot i] \to X$ with the following properties:

\noindent (i)if $a,b$ are such that $f|_{[a,b]}=\operatorname{const.}$ then 
$\gamma_i|_{[a\cdot i,b\cdot i]}$
is a geodesic segment;

\noindent (ii) if $a$ is a point of discontinuity of the function $f$ 
then the  angle between 
$[O,\gamma_i(a\cdot i)]$ and $[\gamma_i(a\cdot i),\gamma_i((a+\epsilon)i)]$ 
is equal to $1/100 + f(a)$, where $\epsilon$ is such 
that $f|_{[a,a+\epsilon]}=\operatorname{const};$

\noindent (iii) if $f|_{[0,a]} \equiv f(0)$  then $\gamma(a)\in l_{2\pi f(0)}$.
\subsection*{Step 2.}
There is always a unique way to start constructing $\gamma$ 
due to the property (iii), 
since for every  function there exists a non--zero interval
$[0,\epsilon]$ where it is constant and is equal to $f(0)$. 
If $\gamma$ is constructed on
$[0,r\cdot i)$ it may be continued in a unique way to the point $r$ by 
setting $\gamma_i(ri)=\lim_{z\to ri-0}\gamma_i(z)$.
If $\gamma_i$ is constructed on $[0,r\cdot i]$ for some
$r<\rho$ we may take such $\epsilon>0$
that $f|_{[r,r+\epsilon]}=\operatorname{const.}$ 
There are two possible cases:
if $r$ is a point of discontinuity of $f$ then we apply (ii) in order to
continue $\gamma_i$ up to the point $(r+\epsilon)i$,
and if $f$ is continuous at $r$ we may apply (i).
Therefore it is possible to construct a unique path $\gamma_i$ on the whole 
$[0,\rho i]$. This path is a broken geodesic line with possibly infinite
number of links.

We set  $\{x_i^f\}=\gamma_i(\rho \cdot i)$. 
Let us prove that the correspondence $f\in A \to \{x_i^f\}\subset X$
satisfies Definition \ref{embed}.
\subsection*{Step 3.}
Our next aim is to show that for any $f,g\in A$ we have
$$\frac{d_X(x_i^f,x_i^g)}{i}\to d_A(f,g).$$

Let us prove first that 
\begin{equation}
\label{sim}
\frac{d_X(x_i^f,O)}{i} \to \rho.
\end{equation}
Denote a function $f_r=f|_{[0,r]}$. 
Consider the set $$R=\{r\in [0,\rho)|\, \forall r'\le r 
\quad \frac{d_X(x_i^{f_r},O)}{i} \to r'\}.$$

We want to verify that $R=[0,\rho]$. We do it in the same
way as we have shown the existence of the path $\gamma_i$ in the previous
Step.

Since for some  $\epsilon>0$ $f|_{[0,\epsilon]}=\operatorname{const.}$
therefore $\epsilon \in R$.
If every  $r_1 <r$ belongs to $R$ then $r\in R$.
This follows from the triangle inequality:
$$
r \ge \frac{d_X(x_i^{f_r}, O)}{i} \ge \frac{d_X(x_i^{f_{r_1}}, O)}{i}- (r-r_1),
$$
and since we may choose $r_1$ arbitrary close to $r$ we get
$$\frac{d_X(x_i^{f_r}, O)}{i}\to r.$$
Finally, if $r \in R$ for some $r<\rho$ then $r+\epsilon \in R$
for some $\epsilon>0$.
Indeed, let $r'$ be the left end of the interval containing $r$ 
on which $f$   is constant.
Due to (ii) the angle $\angle ([O, \gamma_i(r')], 
[\gamma_i(r'), \gamma_i(r+\epsilon)]) = 1/100 + f(r') > 0$ and hence
using Lemma \ref{nep} we obtain that  
$d_X(x_i^{f_{r+\epsilon}}, O)/i \to r+\epsilon$.
This completes the proof of the fact that $R=[0,\rho]$.
\subsection*{Step 4.}
In order to proceed we need the following lemma.
\begin{lemma}
\label{29}
Let $A_iB_iC_i$ be a triangle in $X$ such that $|A_iB_i|=u\cdot i$,
$|A_iC_i|=v \cdot i$, $u,v >0$.
Consider a path $\zeta_i$ of length $|\zeta_i|$ 
connecting $C_i$ and $B_i$ such that
\begin{equation}
\label{cond}
\lim_{i\to \infty}\frac{u\cdot i}{|\zeta_i|+v\cdot i}=1.
\end{equation}
Then $\lim_{i\to\infty}\angle A_i = 0$.
\end{lemma}
\noindent {\bf Proof.} Suppose that the angle $\angle A_i$ does not
tend to zero. Then there exist a subsequence $i_k\to \infty$ such
that $\angle A_{i_k}>\alpha>0$.
Then due to Lemma \ref{nep} $d_X(B_{i_k},C_{i_k})/i_k \to u+v$. 
On the other hand, $|\zeta_{i_k}|\ge d_X(B_{i_k},C_{i_k})$ and hence
$$\lim_{i_k\to \infty}\frac{|\zeta_{i_k}|+v\cdot i_k}{u\cdot i_k} \ge
\frac{u+2v}{u} >1, $$
which contradicts with (\ref{cond}). This completes the proof of 
Lemma \ref{29}. \qed
\subsection*{Step 5.}
Consider a function
$f\in A$ which is discontinious at the point $s$.
We say that it has a 
{\it discontinuity of the first type} at $s$
if there exists $\delta>0$ such that $f_{[s-\delta,s]} \equiv 
\operatorname{const.}$ 
Otherwise we say that $f$ has a {\it discontinuity of the second type} at $s$.
Actually, discontinuities of the second type are accumulation points
of discontinuities of the first type.

Now, let $f,g \in A$ be arbitrary functions and $s$ be their separation
moment. At least one of the functions 
should be discontinuous at $s$, let it be the function $f$.
Note that if it has a discontinuity 
of the first type at $s$, then $g$ is
either continuous, or has a discontinuity of the first type, and
if $f$ has a discontinuity of the second type at $s$, the function $g$
has also a discontinuity of the second type at this point.

One may check in each case that due to the property (ii) and 
Lemma \ref{29} the paths corresponding to functions $f$ and $g$ for large $i$  
separate exactly at the point $x_i^s=x_i^{f_s}=x_i^{g_s}$, and
that the angle
$\angle([x_i^s,x_i^f],[x_i^s,x_i^g)$
is separated from zero. Indeed, if $g$ is continuous at $s$ then $f$ has
a discontinuity of the first type (see Figure 1, left); the angles marked with
black are close to zero by Lemma \ref{29}, 
the angle $\alpha$ by property (ii) is neither $0$ nor $\pi$  and does not
depend on $i$. 
\begin{figure}
\noindent\centering\hfill\epsffile{fig1.1}\hfill\epsffile{fig2.1}\hfill
\caption{}
\end{figure}

If $g$ is discontinuous at $s$ (see Figure 1, right) then similarly 
the angles marked with black
are  close to zero by Lemma \ref{29} and the angle $\beta$ by property (ii)
is neither $0$ nor
$\pi$ and does not depend on $i$.
Thus it can be shown analogously to (\ref{sim}) that
$d_X(O,x_i^s)/i \to~s$, $d_X(x_i^s,x_i^f)/i \to \rho_f-s$,
$d_X(x_i^s,x_i^g)/i \to \rho_g-s$.
Therefore, applying Lemma \ref{nep} we finally obtain:
$$
\lim_{i\to \infty}\frac{d_X(x_i^f,x_i^g)}{i}=
\lim_{i\to \infty}\frac{d_X(x_i^f,x_i^s)+d_X(x_i^s,x_i^g)}{i}=
(\rho_f-s)+(\rho_g-s)=d_A(f,g),
$$
which completes the proof of the theorem. \qed

\subsection* {Acknowledgments.} We would like to thank
Professor A.~Shnirelman for fruitful discussions which stimulated us
to write this paper. We are also grateful to  Professors B.~Bowditch, 
F.~Paulin, L.~Polterovich and O.~Schramm for helpful remarks,
and to V.~Kondrat'ev for TeX-ing Figure 1.
\section*{References}

\noindent [Bea] A. Beardon, The geometry of discrete groups, Springer--Verlag, 
1983

\noindent [Ber] V.N. Berestovskii, Quasicones at infinity
of Lobavhevski spaces, Announcement on the Maltsev 
International  Algebraic Conference, Novosibirsk, 1989.

\noindent [Bes] M. Bestvina, ${\Bbb R}$-trees in topology, geometry, and group
theory, Preprint, 1999.

\noindent [DrW] L. van den Dries, A.J. Wilkie, On Gromov's theorem concerning
groups of polynomial growth and elementary logic, J. of Algebra 89,
1984, 349-374.

\noindent [Dru] C. Drutu, R\'eseaux des groupes semisimple et 
invariants de quasi-isometrie, Ph.D. Thesis, Universit\'e d'Orsay, 1996.

\noindent [DP] A.G. Dyubina, I.V. Polterovich,
Structures at infinity of hyperbolic spaces, Russian Mathematical
Surveys, vol. 53, No. 5, 1998, 239-240.

\noindent [GhH] E. Ghys, P. de la Harpe, Sur le groupes hyperboliques apres
Mikhael Gromov, Birkh\"auser, 1990. 

\noindent [Gr1] M. Gromov, Hyperbolic Groups, in: Essays in group theory,
ed. S.M.Gersten, M.S.R.I. Publ. 8 , Springer--Verlag, 1987, 75-263.

\noindent [Gr2] M. Gromov, Asymptotic invariants of infinite groups, 
Geometric group theory. Vol. 2 (Sussex, 1991), London Math. Soc. Lecture
Note Ser. 182, Cambridge Univ. Press, 1993, 1-295.

\noindent  [Kap] I. Kapovich, A non-quasiconvexity embedding theorem for
hyperbolic groups, Math. Proc. Camb. Phil. Soc. 127 (1999), no. 3, 
461-486. 

\noindent [KapL] M. Kapovich, B. Leeb, On asymptotic cones and quasi-isometry
classes of fundamental groups of 3-manifolds, GAFA, vol. 3, No. 5,
1995, 583-603.

\noindent [MNO] J. Mayer, J. Nikiel, L. Oversteegen, Universal
spaces for ${\Bbb R}$-trees, Trans. AMS, vol. 334, No. 1, 1992, 411-432.

\noindent [N]  J. Nikiel, Topologies on pseudo-trees and applications,
Memoirs AMS, No. 416, 1989. 

\noindent [PSh] I. Polterovich, A. Shnirelman, An asymptotic subcone       
of the Lobachevskii plane as a space of functions, Russian Math. Surveys,  
vol. 52 No. 4, 1997, 842-843.                                                

\noindent [Sh] A. Shnirelman, On the structure of asymptotic space of the 
Lobachevsky plane, Preprint, 1997, 1-23.

\noindent [ThV] S. Thomas, B. Velickovic, Asymptotic cones of finitely 
generated groups, Bull. London. Math. Soc., vol. 32, no. 2, 2000, 203-208. 
\end{document}